\documentclass[english]{amsart}
\usepackage{amscd,amsfonts,amssymb,amsmath}

\newcommand{\Rm}{\mathbb{R}}

\newcommand{\mC}{\ensuremath{\mathcal{C}}}
\newcommand{\mG}{\ensuremath{\mathcal{G}}}

\newcommand{\mY}{\ensuremath{\mathcal{Y}}}
\newcommand{\mH}{\ensuremath{\mathcal{H}}}

\newcommand{\mT}{\ensuremath{\mathcal{T}}}
\newcommand{\mS}{\ensuremath{\mathcal{S}}}
\newcommand{\mF}{\ensuremath{\mathcal{F}}}
\newcommand{\Nm}{\ensuremath{\mathbb{N}}}

\newcommand{\mA}{\ensuremath{\mathcal{A}}}
\newcommand{\mP}{\ensuremath{\mathcal{P}}}
 
\newcommand{\mB}{\ensuremath{\mathcal{B}}}
\newcommand{\Tm}{\ensuremath{\mathbb{T}}}

\newtheorem{lem}{Lemma}
\newtheorem{thm}[lem]{Theorem}

\newtheorem{prop}[lem]{Proposition}
\newtheorem{defn}[lem]{Definition}

\def\qed {\mbox{}\hfill {\small \fbox{}} \\}  
\def\lto{\longrightarrow}
\def\lmto{\longmapsto}

\def\leq{\leqslant}
\def\geq{\geqslant}

\title
{Young measures, superposition and transport}

\author{Patrick Bernard}
\address{
{\rm Patrick Bernard}\\
CEREMADE\\
Universit\'e de Paris Dauphine\\
Pl. du Mar\'echal de Lattre de Tassigny\\
75775 Paris Cedex 16\\
France}
\email{Patrick.Bernard@ceremade.dauphine.fr}

\date{November 2006}
\begin{document}
\begin{abstract}
We discuss a space of Young measures in connection
with some variational problems. We use it to present 
a proof of the Theorem of Tonelli on the existence of
minimizing curves.
We generalize a recent result of Ambrosio,
Gigli and Savar\'e on the decomposition of the weak solutions
of the transport equation.
We also prove, in the context of Mather theory,
the equality between Closed measures and Holonomic measures.
\end{abstract}

\maketitle

Published in: Indiana University Math Journal, 57 No 1 (2008) pp 247--276.
\section{Introduction}
It is by now a well understood  fact that Young measures 
are a very useful tool in variational problems.
In his book, Young exposes how the relaxation to
appropriate spaces of Young measure allow to treat
with great elegance
the problem of lengh-minimizing curves.
In  the present paper, we present 
an extension of Young's approach to the non-parametric situation,
and describe some applications.
This provides a new proof of the theorem of Tonelli
on the existence of curves minimizing a fiberwise convex action.
The objects which appear in this program are related to
some dynamical optimal transportation problems and
to a variational approach of the Euler equation due to 
Arnold and Brenier, see
\cite{BB,BeBu,Br:03,Br:89,Br:93,Br:99}.
Our initial motivation has been to clarify our understanding
of these objects.

We expose in section \ref{s2} and \ref{s3}
the definition and main properties of the measures
 we will work with:
Young measures, transport measures and generalized curves.

It should come as a reward and as an indication of the usefulness
of this theory that we can provide in section \ref{tonelli}
a short and, we believe, elegant
 proof of the famous theorem of Tonelli on the existence
of action-minimizing curves.
We also underline the formal similarity between 
the problem of action minimizing curves and some dynamic optimal
transportation problem as discussed in \cite{BeBu} and other papers.
We obtain  general existence results for these questions.
In order to study the minimizing measures in a general framework,
we need to pursue the study of transport measures.

This is what we do in section  \ref{s5}, 
where we 
state, discuss and prove Theorem \ref{decomposition},
which is certainly the most important result of the present paper.
We call it Young's superposition principle for it is
directly  inspired
by a  result which appears in the appendix of Young's book.
We propose some applications to the continuity equation 
and to the decomposition of optimal transport measures,
that is to the full understanding of the relation between
dynamic optimal transportations and action-minimizing curves. 
It should be noted that although Young's superposition principle
is more general than another  superposition principle
recently obtained by  Ambrosio, Gigli and Savar\'e in
\cite{AGS} many  of its application to the study of 
transport measures minimizing the action defined by a fiberwise
convex integrand
could in fact be obtained from this
especially important  particular case.

In section \ref{s6}, we apply and adapt the ideas of Theorem
\ref{decomposition} to study the closed measures which 
appear in  Mather's theory of minimizing measures.
In \cite{Ma:91} Mather introduced and studied invariant measures
of a Lagrangian system which minimize the action.
These measures turn out to have remarkable property. 
Later, Ma\~n\'e introduced a class of probability measures,
Holonomic measures,
which contain the invariant measures of all 
Lagrangian flows, and  which have the property that 
minimizing closed measures are invariant.
Then Bangert introduced the larger class of closed measure and 
proved, for some specific Lagrangians, that 
minimizing closed measures are invariant. This was generalized
by Fathi and Siconolfi to a much larger class of $C^2$
Lagrangians.
Young's superposition principle allows to generalize these results to
non-regular integrands (with the appropriate definition
of invariance). We also  prove that the holonomic measure
of Ma\~n\'e and the closed measures of Bangert are the same
objects. 
We finish with some generalities of measure theory in the appendix.

I thank Alessio Figalli and Boris Buffoni for their help
at different stages of the elaboration of the present work.

I finished to write this paper in De Giorgi center, Pisa.
This was an occasion to visit the beautiful Camposanto.
There, in a corner, under a scaffolding,
is the sober grave of Leonida Tonelli, 1885-1946,
Accademia dei Lincei.

\section{Young measures}\label{s2}
We define the space of Young measure we will use,
and recall some general results on the topology of this space.
Let $(X,d)$ be a complete and separable metric space.
We denote by $(\mP_1(X),d)$
the Kantorovich-Rubinstein space of 
Borel probability measures on $X$ with finite first moment,
see the appendix.
 Recall that $(\mP_1(X),d)$ is a complete and separable 
metric space.
Let $I=[a,b]$ be a compact interval and let  $\lambda$ be
the normalized Lebesgue measure on $I$.
We denote by 
$\tilde \mY_1(I,X)$ the set of measurable maps
$$
I\ni t\lmto \mu_t\in \mP_1(X).
$$
There is a natural map
$$
\eta_t
\lmto \lambda\otimes \eta_t
$$
from $\tilde \mY_1(I,X)$
to $\mP_1(I\times X)$,
where we denote by $\lambda\otimes \eta_t$
the only measure which satisfies
$$
\int_{I\times X}f(t,x) d(\lambda\otimes \eta_t)(t,x)
=
\int_I \int_X f(t,x) d\eta_t(x) d\lambda (t)
$$
for each bounded Borel function $f:I\times X\lto \Rm$.
The disintegration theorem states that the image of this
map is the set $\mY_1(I,X)$
of probability measures $\eta\in \mP_1(I\times X)$
whose marginal on the component $I$ is the measure $\lambda$.
We call these measures Young measures.
Moreover, two elements 
of $\tilde \mY_1(I,X)$ have the same image if and only if they are
almost everywhere equal.
Note that $\mY_1(I,X)$ is a closed subset of the
 Kantorovich-Rubinstein space $\mP_1(I\times X)$. 
 We endow it from now on
with the induced distance.
The map 
\begin{equation}\label{int}
\eta\lmto \int_{I\times X} f(t,x)d\eta
\end{equation}
is continuous on $\mY_1(I,X)$
for all continuous function $f(t,x):I\times X\lto \Rm$
such that $|f(t,x)|/(1+d(x_0,x))$ is bounded for some
$x_0\in X$.
This continuity holds for many more functions $f$.

\begin{defn}
A Caratheodory integrand is a Borel function
$f(t,x):I\times X\lto \Rm$
which is continuous in the second variable.
A normal integrand is a Borel function 
$f(t,x):I\times X\lto (-\infty,\infty]$
which is lower semi-continuous in the second variable.
\end{defn}

\begin{prop}
The map 
(\ref{int}) is continuous on 
$\mY_1(I,X)$ if $f$ is a Caratheodory integrand such that 
$|f(t,x)|/(1+d(x_0,x))$ is bounded for some $x_0\in X$.
It is lower semi-continuous if $f$ is a normal integrand 
such that $f(t,x)/(1+d(x_0,x))$ is bounded from below.
\end{prop}
\proof
We follow \cite{BeLa:73}, Lemma II.1.1, p 142
for the first part.
By the Scorza-Dragoni Theorem, (see 
\cite{BeLa:73}, Theorem I.1.1, p 132.)
there exists a sequence $J_n$ of compact subsets on $I$
such that $f$ is continuous on $J_n\times X$ and such that 
$\lambda (J_n)\lto 1$ as $n\lto \infty$.
Then, there exists a sequence of continuous functions
$f_n$ such that $|f_n(t,x)|/(1+d(x_0,x))$ is bounded,
independently of $n$, and such that $f_n=f$
on $J_n\times X$.
It follows that the map (\ref{int})
is the uniform limit
of the continuous maps $\eta\lmto \int f_nd\eta$,
and therefore it is continuous.

In order to prove the second part of the statement,
we first write the integrand $f(t,x)=(1+d(x_0,x))g(t,x)$
with a normal integrand $g$ which is bounded from below. 
Then $g$ is  the increasing pointwise limit
of a sequence $g_n$  of bounded Caratheodory integrands,
see \cite{BeLa:73}, Theorem I.1.2, p 138.
Finally,  the map (\ref{int}) is the increasing limit
of the continuous maps $\eta\lmto \int (1+d(x_0,x))g_n(t,x) d\eta(t,x)$, and therefore it 
is lower semi-continuous.
\qed

\begin{thm}\label{compact}
Let $f(t,x)$ be a normal integrand.
Assume that there exists a proper function $l:X\lto [0,\infty)$
and an integrable function $g:I\lto \Rm$
such that $f(t,x)\geq l(x)(1+d(x,x_0))+g(t)$.
Then for each $C\in \Rm$ the set of Young measures 
$\eta\in \mY_1(I,X)$ which satisfy $\int fd\eta \leq C$
is compact.
\end{thm}
\proof
Since the map $\eta\lmto \int fd\eta$
is lower semi-continuous, it is enough to prove that 
the set of Young measures $\eta$
which satisfy 
$\int l(x)d(x,x_0)d\eta \leq C$ is compact.
This set is obviouly $1$-tight, see the Appendix.
\qed

\section{Transport measures and generalized curves}\label{s3}
In the present section, we set $X=TM$, where $M$ is a complete 
Riemannian manifold without boundary.
We endow this tangent space $TM$ with a complete distance $d$
such that the quotient 
$$
\frac{1+d((x_0,0),(x,v))}{1+\|v\|_x}
$$
and its inverse are 
 bounded on $TM$  for one (and then any) point $x_0\in M$.
 The discussions below do not depend on the choice of this 
 distance $d$.
In order to prove that such a distance exists,
we can  isometrically embed $M$
into a Euclidean space $\Rm^d$  and restrict the distance
$$
D((x,v),(x',v'))=
\min(1,|x'-x|)+|v'-v|,
$$
where $|.|$ is the Euclidean norm on $\Rm^d$.
We fix a compact interval $I=[a,b]$.
and  denote by $C_1(I\times TM)$ the set of continuous functions
$f:I\times TM\lto \Rm$ such that 
$$
\|f\|_1:=\sup_{(t,x,v)\in I\times TM} \frac {\|f(t,x,v)\|}{1+\|v\|_x}<\infty.
$$

\begin{defn}
A \emph{transport measure} is a measure $\eta \in \mY_1(I,TM)$
which satisfies the relation
\begin{equation}\label{transporteq}
\int_{I\times TM} \partial_t g+\partial_xg\cdot v \;d\eta(t,x,v)=0
\end{equation}
for all smooth compactly supported functions 
$g:]a,b[\times M\lto \Rm$.
We denote by $\mT(I,M)\subset \mY_1(I,TM)$ the set of all transport measures. 
Given two probability measures  $\mu_i$ and $\mu_f$
on $M$,
we say that the transport measure $\eta$ is a transport measure
between $\mu_i$ and $\mu_f$
if, in addition, we have 
$$
\int_{I\times TM} \partial_t g+\partial_xg\cdot v \;d\eta(t,x,v)=
\int_{M} g_b(x) d\mu_f(x)-\int_{M} g_a(x) d\mu_i(x)
$$
for each smooth compactly supported 
 function $g:[a,b]\times M\lto \Rm$.
 We denote by $\mT_{\mu_i}^{\mu_f}(I,M)$ the set of transport measures
 between $\mu_i$ and $\mu_f$.
\end{defn}

Note that $\mT(I,M)$ and $\mT_{\mu_i}^{\mu_f}(I,M)$ 
are closed subsets of $\mY_1(I,M)$.
Recalling that we denote by $\mP(M)$ the set of Borel probability
measures endowed with the narrow topology, we have :

\begin{lem}
Let $\eta\in \mT(I,M)$ be a 
transport measure. There exists a continuous family
$\mu_t:I\lto \mP(M)$ of probability measures on $M$
and a disintegration $\eta_t\in \tilde \mY_1(I,TM)$ of $\eta$
such that, for each $t$,  $\mu_t$ is the marginal of $\eta_t$
on the base $M$. We then have $\eta\in \mT_{\mu_a}^{\mu_b}(I,M)$.
\end{lem}

\proof
Let us choose a disintegration $\eta_t$ of $\eta$,
and let $\tilde \mu_t$ be the marginal of $\eta_t$ on $M$.
We want to prove that there is a narrowly  continuous 
map $\mu_t:I\lto \mP(M)$
which is equal to $\tilde \mu_t$ for almost each $t$.
In view of general remarks recalled in the Appendix, it is enough
to prove that, for each smooth and compactly supported function
$f:M\lto \Rm$, the function $t\lmto F(t):=\int fd\mu_t$
is equal almost everywhere to a continuous function.
By applying the equation (\ref{transporteq}) 
to functions $g(t,x)=\phi(t)f(x)$, we get  that 
$F'(t)=\int df_x\cdot v d\eta_t(x,v)$
in the sense of distributions.
It imples that the function $F$ is equal almost everywhere
to an absolutely continuous function.
\qed

\begin{lem}
Let $g(t,x):I\times M\lmto \Rm$
be a $C^1$  bounded and Lipschitz function.
Then for each interval $[\alpha,\beta]\subset [a,b]$,
we have
\begin{equation}\label{transeq2}
\int_{[\alpha,\beta]\times TM} \partial_tg+\partial_xg\cdot v d\eta
=
\int_M g_{\beta}d\mu_{\beta}-
\int_M g_{\alpha}d\mu_{\alpha}
\end{equation}
\end{lem}

\proof
Let us first assume that $g$ is a smooth compactly supported function.
Let us set $F(t)=\int g_td\mu_t$.
It is easy to prove using 
(\ref{transporteq}) that 
$$
F'(t)=
\int_{TM} \partial_t g+ \partial_x g\cdot v d\eta_t
$$
in the sense of distribution.
The desired equality follows by integration.
If $g$ is $C^1$ and compactly supported, then we prove 
(\ref{transeq2})
by approximating $g$ by smooth compactly supported functions.
Let us expose a bit more carefully how the equality can be extended
to bounded and Lipschitz functions which
are not necessarily compactly supported.
We consider an increasing  sequence $\xi_n: M\lto [0,1]$ of
 smooth equi-Lipschitz 
compactly supported functions
such that, for each relatively compact  open set $U$,
we have $\xi_n=1$ on $U$ after a certain rank.
Then (\ref{transeq2}) holds 
for the function $g\xi_n$:
$$
\int_{[\alpha,\beta]\times TM} 
\xi_n\partial_tg+\xi_n\partial_xg\cdot v 
+g\partial_x \xi_n\cdot v d\eta
=
\int_M g_{\beta}\xi_nd\mu_{\beta}-
\int_M g_{\alpha}\xi_nd\mu_{\alpha}.
$$
Thanks to the dominated convergence theorem, we get 
(\ref{transeq2}) at the limit.
\qed

\begin{defn}
The transport measure  $\eta $ is called a generalized curve 
if $\mu_t$ is a dirac measure for each $t\in I$.
Then, there exists a continuous curve $\gamma(t):I\lto M$
such that $\mu_t=\delta_{\gamma(t)}$ for each $t$.
We say that $\eta$ is a generalized curve above $\gamma$.
We denote by $\mG(I,M)$ the set of generalized curves.
\end{defn}

A continuous curve $\gamma:I\lto M$ is  absolutely
continuous if and only if the function $\phi\circ \gamma:I\lto \Rm$
is absolutely continuous for each smooth 
and compactly supported function $\phi:M\lto \Rm$.
We denote by $W^{1,1}(I,M)$ the set of absolutely continuous
 curves.
 We say that a sequence $\gamma_n$ is converging to $\gamma$
 in $W^{1,1}(I,M)$ if the function sequence 
 $$
 d\big((\gamma_n(t),\dot \gamma_n(t)), (\gamma(t),\dot \gamma(t))
 \big)
 $$
 is converging to zero in $L^1$, or equivalentely if
 the following three conditions are satisfied:
 \begin{itemize}
 \item
The sequence  $\gamma_n$ is converging uniformly to $\gamma$.
\item
The sequence $(\gamma_n(t),\dot \gamma_n(t)):I\lto TM$ is converging in measure to
$(\gamma(t),\dot \gamma(t))$.
\item
The sequence $\|\dot \gamma_n(t)\|_{\gamma_n(t)}$ is equi-integrable, 
or equivalentely it is relatively weakly  compact in $L^1(I,\Rm)$.
\end{itemize}
It is well-known that smooth curves are dense in $W^{1,1}(I,M)$.
\begin{lem}
Let $\Gamma\in \mT(I,M)$ be a generalized curve.
Then there exists an absolutely continuous curve $\gamma(t)$
such that $\Gamma$ is a generalized curve above $\gamma$ and
 there exists a measurable family $\Gamma_t$ of 
probabilty measures
on $T_{\gamma(t)}M$ such that 
 $\Gamma=dt\otimes \delta_{\gamma(t)}\otimes \Gamma_t$,
 which means that
 $$
 \int_{I\times TM} f(t,x,v) d\Gamma(t,x,v)
 =
 \int_I \int_{T_{\gamma(t)}M} f(t,\gamma(t),v)\, d\Gamma_t(v) \,dt
 $$
 for each $f\in L^1(\Gamma)$.
  In order that this formula defines a generalized curve
 above the absolutely continuous curve $\gamma$, it is necessary
 and sufficient that the function 
 $t\lmto \int_{T_{\gamma(t)}M}\|v\|_{\gamma(t)}d\Gamma_t(v)$
 is $\lambda$-integrable on $I$, and that 
 $\int_{T_{\gamma(t)}M}vd\Gamma_t(v)=\dot\gamma(t)$
for almost all $t$.
\end{lem}
\proof
Let $\Gamma$ be a generalized curve over $\gamma$.
By the disintegration theorem, the measure $\Gamma$
can be written on the form 
$\Gamma=dt\otimes \delta_{\gamma(t)}\otimes \Gamma_t$ 
with some measurable family $\Gamma_t$ of probability measures
on $T_{\gamma(t)}M$. 
We want to prove that the curve $\gamma(t)$ is absolutely continuous
and that 
$$
\dot \gamma(t)=\int_{\Rm^d}vd\Gamma_t(v)
$$
for almost all $t$.
It is enough to prove that, for each smooth compactly supported function $\phi:M\lto \Rm$, we have
$$
(\phi\circ \gamma)'(t)=
d\phi_{\gamma(t)}\cdot \int_{T_{\gamma(t)}M}vd\Gamma_t(v)
$$
in the sense of distributions.
For each smooth compactly supported function
$
f(t): ]a,b[\lto \Rm
$ 
we can apply the equation (\ref{transporteq})
to the function $g(t,x)=f(t)\phi(x)$,
and get
\begin{eqnarray*}
0&=&
\int_{I\times TM} f'(t)\phi(x)+f(t)d\phi_x\cdot v d\Gamma(t,x,v)
\\
&=&\int_0^1 f'(t)\phi(\gamma(t)) dt + 
\int_0^1 f(t) \int_{T_{\gamma(t)}M} d\phi_{\gamma(t)}\cdot v\; d\Gamma_t(v)\; dt.
\end{eqnarray*}
This implies that $\phi\circ \gamma$
is absolutely contiuous and that 
$$
(\phi\circ \gamma)'(t)=
 \int_{T_{\gamma(t)}M} d\phi_{\gamma(t)}\cdot v\; d\Gamma_t(v)
 =d\phi_{\gamma(t)}\cdot \int_{T_{\gamma(t)}M}vd\Gamma_t(v)
$$
which is the desired result.
\qed

\begin{thm}\label{closed}
The set  $\mG(I,M)$ of generalized curves is closed  
in $\mY_1(I,TM)$.
In addition, the map $\mG \lto C^0(I,\Rm^d)$
which, to a generalized curve $\Gamma$ above $\gamma$,
associates the curve $\gamma$ is continuous.
\end{thm}

\proof
Let $\Gamma_n$ be a sequence of generalized curves
converging in $\mP_1(I\times TM)$ to a limit $\eta$.
We have to prove that $\eta$ is a generalized curve.
The family $\eta,\Gamma_1,\Gamma_2, \ldots, \Gamma_n,
\ldots$ is compact in $\mP_1(I\times TM)$, hence it 
has uniformly integrable first moment.
This implies that the sequence $\gamma_n$ of associated
curves is equi-absolutely continuous.
Taking a subsequence, we can assume that the sequence 
$\gamma_n$
has a limit $\gamma$ in $C^0(I,M)$.
It is not hard to check, then, that 
$\eta$ is a generalized curve 
above $\gamma$.
\qed

If $\gamma:I\lto M$  is absolutely continuous,
then we will denote by $\bar \gamma$
the generalized curve above $\gamma$
given by
$$
\int_X f(t,x,v) d\bar \gamma(t,x,v)
=\int_I f(t,\gamma(t),\dot \gamma(t)) dt
$$
for each bounded Borel function $f$.
In other words, we have 
$$\bar \gamma=dt\otimes \delta_{\gamma(t)}\otimes
\delta_{\dot \gamma(t)}.
$$
We denote by $\mC(I,M)\subset \mT(I,M)$ the set of  transport measures which
are of that form. 

\begin{lem}\label{continu}
The map 
$$
W^{1,1}(I,M)\lto \mG(I,M)
$$
$$
\gamma\lmto \bar \gamma
$$
is continuous.
\end{lem}

\proof
Let  $\gamma_n\in W^{1,1}(I,M)$
be a sequence which converges to $\gamma$.
We have to prove that 
$$
\int_I f(t,\gamma_n(t),\dot \gamma_n(t)) d\lambda 
\lto \int_I f(t,\gamma(t),\dot \gamma(t)) d\lambda 
$$
for each $f\in C_1(I\times TM)$.
Since the sequence 
$(\gamma_n(t),\dot \gamma_n(t))$ is converging in measure 
to $(\gamma(t), \dot \gamma(t))$, 
we can suppose by extracting a subsequence
that it is converging almost everywhere.
The desired convergence follows from the observation
that the sequence of real functions
$$
t\lmto f(t,\gamma_n(t),\dot \gamma_n(t)) 
$$
is converging almost everywhere to 
 $f(t,\gamma(t),\dot \gamma(t))$
 and is equi-integrable because 
 $$
 |f(t,\gamma_n(t),\dot \gamma_n(t))  |
 \leq \|f\|_1 \big(1+\|\dot \gamma_n(t)\|_{\gamma_n(t)}\big)
  $$
  and, by definition of the convergence in $W^{1,1}$,
  the sequence $\|\dot \gamma_n(t)\|_{\gamma_n(t)}$
  is equi-integrable.
 \qed

Let us mention  for completeness:

\begin{thm}
The set $\mG(I,M)$ of generalized curves is the closure,
in $\mY_1(I,TM)$, of the set $\mC(I,M)$ of curves.
\end{thm}

\section{Tonelli Theorem and optimal transportation}\label{tonelli}
In the present section,
we use transport measures and generalized curves to expose 
some  results on the existence of certain minimizers.
The results are well-known, but the presentation is somewhat original.
 We consider 
a normal integrand $L:[a,b]\times TM\lto \Rm \cup \{+\infty\}$.
We say that $L$ is fiberwise convex if, for each 
fixed $(t,x)$, the function $v\lmto L(t,x,v)$
is convex on $T_xM$.
The role of convexity in minimization problems is enlightened by
the following standard observation:

\begin{lem}\label{cx}
Let $L$ be a fiberwise convex normal integrand.
If $\Gamma$ is a generalized curve above $\gamma$,
then
$$
\int Ld\Gamma\geq \int_a^b L(t,\gamma(t),\dot \gamma(t))dt
=\int Ld\bar \gamma.
$$
\end{lem}

\proof
For each $t$, we have 
$$
\int_{T_{\gamma(t)}M} L(t,\gamma(t),v)d\Gamma_t(v)
\geq L(t,\gamma(t),\dot \gamma(t))
$$
by Jensen inequality.
We obtain
$$
\int Ld\Gamma=
\int_0^1 \int_{T_{\gamma(t)}M} L(t,\gamma(t),v)d\Gamma_t(v)dt
\geq \int_0^1 L(t,\gamma(t),\dot \gamma(t))dt=\int Ld\bar \gamma.
$$
\qed
We now discuss the classical problem of the existence of minimizing curves.
We fix two points
 $x_i$ and $x_f$ in  $M$, and
consider the set  $AC_{x_i}^{x_f}$
of absolutely continuous curves 
$\gamma:I\lto M$ such that $\gamma(a)=x_i$ and $\gamma(b)=x_f$.
We also consider the set 
$$
\mG_{x_i}^{x_f}=\mG(I,M)\cap \mT_{\delta_{x_i}}^{\delta_{x_f}}(I,M)
$$
of generalized curves above elements of $AC_{x_i}^{x_f}$.
Note that $\mG_{x_i}^{x_f}$ is closed in 
$\mT(I,M)$.
The action of an absolutely continuous curve $\gamma$
is the integral
$$\int_a^b L(t,\gamma(t),\dot \gamma(t))dt,
$$
the action of a transport measure $\eta$
is the integral $\int_{I\times TM} Ld\eta$.
The following result is well-known:

\begin{thm}\label{T0}
Let $L(t,x,v):[a,b]\times TM\lto \Rm\cup {+\infty}$ be a normal integrand. 
We assume that the integrand $L$ satisfies:
\begin{itemize}
\item[(L1)]
the quotient
$$
\frac{L(t,x,v)}{1+\|v\|_x}
$$
is bounded from below and proper.
\end{itemize}
For each $C\in \Rm$ the set
$$
\mA_C^g:=\{\Gamma\in \mG_{x_i}^{x_f} 
\quad |\quad \int Ld\Gamma\leq C\}\subset \mG_{x_i}^{x_f}
$$
is compact, and if $L$ is fiberwise convex,  the set
$$
\mA_C:=\{\gamma\in AC_{x_i}^{x_f} \quad |\quad 
\int_a^bL(t,\gamma(t),\dot \gamma(t))\leq C\}\subset C(I,M)
$$
is compact for the uniform topology.
\end{thm}

As a major consequence, we obtain that 
the action reaches its minimum on  $\mG_{x_i}^{x_f}$
if there exists a generalized curve of finite action 
in  $\mG_{x_i}^{x_f}$.
If in addition the integrand is fiberwise convex, 
then the action also reaches its minimum
on $AC_{x_i}^{x_f}$, and we have
$$
\min_{\Gamma\in \mG_{x_i}^{x_f}} \int Ld\gamma
=
\min_{\gamma\in AC_{x_i}^{x_f}} 
\int_a^b L(t,\gamma(t),\dot \gamma(t))dt.
$$

\proof
The compactness of $\mA_C^g$ follows from
Theorem \ref{compact}.
If $L$ is fiberwise convex, then, by Lemma \ref{cx},
the set $\mA_C$ is the image of the compact set
$\mA^g_C$ by the continuous map $\Gamma\lmto \gamma$
(the map which, to a generalized curves $\Gamma$
above $\gamma$, associates the curve $\gamma$).
\qed

In applications, it is useful to have the following stronger and
still standard result:

\begin{thm}[Tonelli]
The same conclusions (as Theorem \ref{T0})
hold if the hypothesis (L1) on the integrand is replaced
by the two following ones:
\begin{itemize}
\item[(L2)]
The integrand $L$ is uniformly superlinear over each compact subset 
of $M$. It means that, for each compact $K\subset M$, there
exists a function 
$l:\Rm^+\lto \Rm$ such that $\lim_{r\lto \infty} l(r)/r=\infty$
and such that $L(t,x,v)\geq l(\|v\|_x)$ for each $(t,x,v)\in [a,b]
\times T_KM$.
\item[(L3)]
There exists a positve constant $c$ such that 
$L(t,x,v)\geq c(\|v\|_x-1)$.
\end{itemize}
\end{thm}

\proof
We have to prove that the set of generalized curves 
$\Gamma \in \mG_{x_i}^{x_f}$ which satisfy 
$\int L d\Gamma\leq C$  is compact.
Using (L3), we see that, if $\Gamma$ is a generalized curve
over $\gamma$, then
$$
\int_a^b \|\dot \gamma(t)\|_{\gamma(t)}dt\leq (C+(b-a))/c.
$$
Let $K$ be the closed ball (for the Riemaniann distance on $M$)
of center $x_i$ and radius $(C+b-a)/c$. This ball is compact
because $M$ is complete.
Let us define a modified integrand $L_K$ by
 $L_K(t,x,v)=L(t,x,v)$ if $x\in K$
and $L_K(t,x,v)=+\infty$ if $x\not \in K$.
 A generalized curves
$\Gamma \in \mG_{x_i}^{x_f}$  satisfy 
$\int L d\Gamma\leq C$
if and only if it satisfies 
$\int L_Kd\Gamma\leq C$.
Since $L_K$ satisfies (L1), we conclude by Theorem \ref{T0}.
\qed

We can extend these
considerations to more general boundary conditions.
Our presentation allows to see the following dynamic optimal transportation problem as a natural generalisation of 
Tonelli theorem.

\begin{thm}
Let $L$ be a normal integrand which satisfies :
\begin{itemize}
\item[(L4)]
The integrand $L$ is uniformly superlinear :
there exists a function 
$l:\Rm^+\lto \Rm$ such that $\lim_{r\lto \infty} l(r)/r=\infty$
and such that $L(t,x,v)\geq l(\|v\|_x)$ for each $(t,x,v)\in [a,b]
\times TM$.
\end{itemize}
Let $\mu_i$ and $\mu_f$ be two Borel probability measures on
$M$. 
Then for each $C\in \Rm$, the set $\mB_C$ of transport measures
 $\eta \in \mT_{\mu_i}^{\mu_f}$
 which satisfy $\int Ld\eta\leq C$
 is compact. 
\end{thm}

Note that (L4) implies (L2) and (L3).

\proof
The conclusion would be obvious if $L$ satisfied (L1),
but (L4) is weaker.
For each $\epsilon>0$, there exists a constant $R$ such that 
$$
\int_{\|v\|_x\geq R} (1+\|v\|_x)d\eta(t,x,v)\leq \epsilon
$$
for each $\eta\in \mB_C$.
This is a direct consequence of (L4).
We claim  that there exists a compact ball $B\subset M$
such that $\eta(I\times T_BM)\geq 1-\epsilon/(1+R)$ 
for each $\eta\in \mB_C$.
Assuming the claim, we have 
$$
\int_{\{(t,x,v)\in I\times TM \;| x\not \in B \text{ or }
 \|v\|_x\geq R\}}
(1+\|v\|_x) d\eta(t,x,v)
\leq 2\epsilon
$$
for each $\eta\in \mB_C$.
Therefore, $\mB_C$ is $1$-tight and thus compact.
Let us now prove the claim.
For each $\Delta>0$, 
there exists a $C^1$, bounded and $1$-Lipschitz function 
$g:M\lto [0,\Delta]$ such that 
$g=\Delta $ outside of a compact ball $B$ and such that 
$\int g d\mu_i \leq 1$.
For $\eta \in \mB_C$, we have 
$$
\int_M g d\mu_t=\int_M g d\mu_i+
\int _{[a,t]\times TM} dg_x\cdot v \,d\eta
\leq 1+(C+b-a)/c
$$
where $c$ the constant of (L3).
We conclude that 
$$
\int g d\eta= \int_a^b \int_M g d\mu_t dt 
\leq (b-a)(1+(C+b-a)/c)
$$
for each $\eta\in \mB_C$. 
It follows that
$$\eta(I\times TM-I\times T_BM)
\leq (b-a)(1+(C+b-a)/c)/\Delta.
$$
Since $\Delta$ can be chosen arbitrarily, the
 claim is proved.
\qed

Some general comments are needed before we can 
describe the additional conclusions satisfied for fiberwise
convex Lagrangians.
If $\eta\in \mY_1(I,TM)$ 
is a Young measure, then we call 
$\mu$ the image of $\eta$ by the projection
 $I\times TM\lto I\times M$.
 We can desintegrate $\eta$ with respect to this projection 
and obtain a measurable family
$\eta_{t,x}$ of probability measures on $T_xM$
such that $\eta=\mu\otimes \eta_{t,x}$.
We define the vectorfield
$V(t,x):I\times M\lto TM$
by the expression
$$
V(t,x)=\int_{T_xM} vd\eta_{t,x}(v).
$$
Note that $V(t,x)$ is a Borel time-dependant vector-field,
and that the integrability condition 
$$\int \|V(t,x)\|_x d\mu(t,x)<\infty
$$
is satisfied.
\begin{lem}
The Young measure $\eta\in \mY_1(I,TM)$
is a transport measure if and only if the continuity equation
\begin{equation}\tag{PDE}\label{PDE}
\partial_t \mu + \text{div}(V\mu)=0.
\end{equation}
holds in the sense of distributions.
\end{lem}

The couple $(V,\mu)$  is what we called in \cite{BeBu}
the transport current asssociated to the transport measure
$\eta$.
Such objects were previously introduced by Benamou and Brenier,
see \cite{BB}, \cite{Br:99} and \cite{Br:03}.

\proof
A test function is a smooth and compactly supported
function on $]a,b[\times M$.
The measure $\eta$ is a transport measure if and only if
$$
\int_{I\times M}\int _{T_xM} \partial_t g(t,x)+\partial_x g(t,x)\cdot v
d\eta_{t,x}(v) d\mu(t,x) =0
$$
for each test function.
The equation (PDE) holds in the sense of distributions 
if and only if
$$
\int_{I\times M} \partial_t g(t,x)+\partial_x g(t,x)\cdot 
V(t,x) d\mu(t,x)=0
$$
for each test function $g$.
The equivalence follows from the observation that
$$
\int_{T_xM} \partial_xg(t,x)\cdot v d\eta_{t,x}
=\partial_xg(t,x)\cdot \int_{T_xM} v d\eta_{t,x}
=\partial_xg(t,x)\cdot V(t,x)
$$
by definition of $V$.
\qed

Conversely,
consider a Borel vector-field $V(t,x):I\times M\lto TM$
and a probability measure $\mu$ on $I\times M$ whose
marginal on $I$ is $\lambda$.
Assume that (PDE) holds and that the integrability condition
$\int \|V\|d\mu<\infty$ is satisfied.
Then, the measure $\tilde V_{\sharp}\mu$ is a transport measure,
where $\tilde V(t,x)=(t,V(t,x))\in I\times TM$.
The following  generalization of Lemma \ref{cx}
is now obvious:

\begin{lem}
Let $L$ be a fiberwise convex normal integrand.
If $\eta$ is a transport measure, and $\mu$ and $V$
are associated to it as above, then 
$\tilde V_{\sharp} \mu$ is a transport measure,
and
$$
\int Ld\eta\geq \int L d(\tilde V_{\sharp} \mu).
$$
As a consequence, if there exists a transport measure minimizing
the action in $\mT_{\mu_i}^{\mu_f}(I,M)$, then there exists 
a minimizing transport measure in $\mT_{\mu_i}^{\mu_f}(I,M)$
which is concentrated on the graph of a Borel
vector-field.
\end{lem}

\section{The superposition principle}\label{s5}
The main stream of this section consists of writing tranport measures
as superpositions of generalized curves.
This is the adaptation to the non-parametric setting
of a theory sketch in the appendix of Young's book.

\subsection{Young's superposition principle.}
We first adapt an important result of Young:

\begin{thm}[Young]\label{Young}
The set $\mT(I,\Rm^d)$ of transport measures is the closed
convex envelop in $\mY_1(I,\Rm^{2d})$ of the set $\mC(I,\Rm^d)$ of curves
(and hence also of the set $\mG(I,\Rm^d)$ of generalized curves).
\end{thm}

Let us immediately 
 mention the restatement  Young's result which we will use:

\begin{thm}\label{decomposition}
If $\eta$ is a transport measure on a complete manifold $M$,
then there exists a Borel measure $\nu$
on  $\mG(I,M)$ such that
$\eta= \int_{\mG} \Gamma d\nu(\Gamma)$, which means that
\begin{equation}\label{S}
\int_{I\times TM} f d\eta= \int _{\mG} 
\int _{I\times TM} fd\Gamma d\nu(\Gamma)
\end{equation}
for each function $f\in L^1(\eta)$. We then say that $\nu$ is a 
decomposition of $\eta$.
\end{thm}

Let us make a few simple remarks before proving these results.

\begin{prop}\label{support}
If $\eta$ is concentrated 
on the Borel subset $Y\subset I\times TM$,
and if $\nu$ is a decomposition of $\eta$,
 then $\nu$-almost every 
generalized curve $\Gamma$ is concentrated on $Y$.
\end{prop}
\proof
Apply
(\ref{S}) with $f=0$ on $Y$ and $f=1$ outside of $Y$.
We get $\int_{I\times TM} fd\Gamma=0$ for $\nu$-almost all 
$\Gamma$, which means that $\Gamma$ is concentrated on 
$Y$.
\qed

For each $t\in I$, let $ev_t :\mG(I,M)\lto M$ be the continuous
map obtained by composing the natural projection
$\mG\lto C^0(I,M)$
and the evaluation map $\gamma\lmto \gamma(t)$.

\begin{prop}
If $\nu$ is a decomposition of $\eta$,
and if $\mu_t$ is the continuous family of probability
measures on $M$ associated to $\eta$, then 
 $\mu_t=(ev_t)_{\sharp}\nu$,
\end{prop}

\proof
We denote by $\gamma_{\Gamma}$ the continuous curve associated
to the generalized curve $\Gamma$.
It is enough to prove that 
$$
\int f(t,x)d\eta= \int_{\mG} f(t, \gamma_{\Gamma}(t)) dt d\nu (\Gamma)
$$
for each continuous and bounded function 
$f:I\times M\lto \Rm$.
This follows from the fact that 
$$
\int_X f(t,x) d\Gamma (t,x,v)=
\int_0^1 f(t,\gamma_{\Gamma}(t))dt 
$$
for each generalized curve $\Gamma$. 
\qed

Finally, let us explain how Theorem \ref{decomposition}
follows from Theorem \ref{Young}.
We isometrically embed the manifold $M$
as a closed subset of some Euclidean space $\Rm^d$.
Then the Transport measures and generalized curves
on $M$ are just the transport measures and generalized curves
on $\Rm^d$ which are supported on 
$I\times TM\subset I\times \Rm^{2d}$.
Let $\eta$ be a transport measure. In view of Young's theorem
and of the appendix, $\eta$ admits a decomposition $\nu$
by generalized curves on $\Rm^d$. 
By Proposition \ref{support} above,
$\nu$ almost every generalized curve $\Gamma$ is supported on 
$I\times TM$, hence $\nu$ can be seen as a probability measure 
on $\mG(I,M)$.

\subsection{Proof of Young's superposition principle}
We prove the superposition principle by duality,
following the sketch of proof proposed by Young in his book.
By Proposition \ref{dual} of the appendix,
it is enough to prove that, for each function 
$f\in C_1(I\times \Rm^{2d})$ such that
\begin{equation}\label{0}
\int_0^1 f(t,\gamma(t),\dot \gamma(t))dt \geq 0 \quad
\forall \gamma\in W^{1,1}(I,\Rm^d)
\end{equation}
we have
$
\int fd\eta\geq 0
$
for all transport measures $\eta\in \mT(I,\Rm^d)$.
It is sufficient to obtain the conclusion for functions
$f\in C_1(I\times TM)$
which satisfy 
\begin{equation}\label{1}
\int_0^1 f(t,\gamma(t),\dot \gamma(t))dt \geq 1\quad
\forall \gamma\in W^{1,1}(I,\Rm^d).
\end{equation}
Indeed, if this is proved, and if $f$ satisfies
(\ref{0}), then for each $\epsilon>0$,
the function $(f+\epsilon)/\epsilon$
satisfies (\ref{1}), hence 
$\int fd\eta \geq -\epsilon$ for each transport measure $\eta$,
and finally $\int fd\eta \geq 0$.

Let us   fix a function $f\in C_1(I\times \Rm^{2d})$,
assume (\ref{1}),
 and define
the value function $u:I\times \Rm^d\lto \Rm$
by
$$
u(t,x):=\inf_{\gamma\in W^{1,1}(\Rm,\Rm^d),\gamma(t)=x}
 \int_0^t f(s,\gamma(s),\dot\gamma(s)) ds.
$$
We have the equality
$$
u(t,x)= \inf_{\gamma\in W^{1,1}(\Rm,\Rm^d),\gamma(t)=x}
u(s,\gamma(s))+
\int_s^tf(\sigma,\gamma(\sigma),\dot\gamma(\sigma)) d\sigma
$$
for each $s\leq t$ and each $x$. This equality is 
called the dynamic programming principle.

\begin{lem}\label{droite}
We have
$$
u(t,y)\leq u(s,x)+\|f\|_1((t-s)+|y-x|)
$$
for each $s\leq t$  in $I$ and each $x,y$ in $\Rm^d$.
\end{lem}
\proof
Just observe that 
$$
u(t,y)\leq u(s,x)+\int_s^t f(\sigma,x+\sigma(y-x)/(t-s),(y-x)/(t-s))d\sigma$$
$$
\leq u(s,x)+(t-s)\|f\|_1(1+|y-x|/(t-s))
$$
\qed

\begin{lem}\label{usc}
The value function $u$ is bounded and upper
semi-continuous. In addition, we have $u(0,x)=0$ and
$u(1,x)\geq 1$ for all $x$.
\end{lem}

\proof
The inequality $u(1,x)\geq 1$ follows from 
(\ref{1}).
For each $ \gamma\in W^{1,1}(\Rm,\Rm^d)$,
 let us consider the function
$$
u_{\gamma}(t,x):= \int_0^t f(s,\gamma(s)+x-\gamma(t),\dot\gamma(s))
ds
$$
which is continuous and bounded.
Observing that $u=\inf_{\gamma\in W^{1,1}(\Rm,\Rm^d) } u_{\gamma}$,
we conclude that the function $u$
is upper semi-continuous and bounded from above.
It follows from Lemma \ref{droite} that 
$u(t,x)\geq u(1,x)+\|f\|_1(t-1)\geq 1+ \|f\|_1(t-1)$ is bounded from below.
\qed

\begin{lem}
There exists sequences $u_n:I\times \Rm^d\lto \Rm$ and 
$f_n:I\times \Rm^{2d}\lto \Rm$ of functions such that:
\begin{itemize}
\item The sequence $f_n$ is bounded in $C_1(I\times \Rm^{2d})$
and $f_n\lto f$ pointwise.
\item The functions $u_n$ are smooth, bounded and Lipschitz.
They satisfy $u_n(0,x)=u_n(1,x)=0$ for all $n$ and all $ x$.
\item
 The inequality
$$
 u_n(t,\gamma(t))-u_n(s,\gamma(s))\leq
\int_s^t f_n(\sigma,\gamma(\sigma),\dot\gamma(\sigma))d\sigma
$$
holds for each $s\leq t$ in $\Rm$ and each absolutely continuous curve
$\gamma:\Rm\lto \Rm^d$.
\end{itemize}
\end{lem}

\proof
There exists $\delta>0$
such that $u(t,x)< 1/2$ when $t\leq a+\delta$
and $u(t,x)> 1/2$ when $t\geq b-\delta$.
It is convenient to consider the function 
$\tilde f_n:\Rm\times \Rm^d\times \Rm^d
\lto \Rm$
which is equal to $f$ on 
$[a+2/n,b-2/n]\times \Rm^d\times \Rm^d$ and to $0$ outside of this
set,
and the function $\tilde u_n:\Rm\times \Rm^d
\lto \Rm$
which is equal to $u-1/2$ on $[a+2/n,b-2/n]\times \Rm^d$,
to $0$ outside of this set.
Note that 
$$
\tilde u_n(t,\gamma(t))-\tilde u_n(s,\gamma(s))\leq
\int_s^t \tilde f_n(\sigma,\gamma(\sigma),\dot\gamma(\sigma))d\sigma
$$
for each $n$, each $s\leq t$ in $\Rm$ and each absolutely continuous curve
$\gamma:\Rm\lto \Rm^d$.

Let $\rho_n(t,x):\Rm\times \Rm^d\lto [0,\infty)$
be a sequence of convolution kernels, that is of smooth 
non-negative functions 
such that $\int_{\Rm\times \Rm^d}\rho_n(t,x)dxdt=1$
and such that $\rho_n$ is supported on the ball of center 
$0$ and radius $1/n$.
Let us define the  functions
$
u_n=\rho_n*\tilde u_n:\Rm\times \Rm^d
\lto \Rm
$:
$$u_n(t,x)=
\int_{\Rm\times \Rm^d} 
\tilde u_n(t-\sigma,x-y)
 \rho_n(\sigma,y)
d\sigma dy,
$$
and
$f_n:\Rm\times \Rm^d\times \Rm^d
\lto \Rm$
by
$$
f_n(t,x,v)=
\int_{\Rm\times \Rm^d}
\tilde f_n(t-\sigma,x-y,v)\rho_n(\sigma,y)
d\sigma dy.
$$
 For each fixed curve $\gamma$ and each $n$, the inequality
$$
\tilde u_n(t-\sigma,\gamma(t-\sigma)-y)-
\tilde u_n(s-\sigma,\gamma(s-\sigma)-y)\leq
\int_{s}^{t} \tilde f_n(\zeta-\sigma,\gamma(\zeta-\sigma)-y,
\dot\gamma(\zeta-\sigma))d\sigma
$$
holds for each $(\sigma,y)$, and then the 
third point of the Lemma is obtained by integration.
\qed

Let $\eta$ be a transport measure.
We want to prove that $\int fd\eta\geq 0$.
Let us set
$$
h_n(t,x,v):=f_n(t,x,v)-\partial_tu_n(t,x)-\partial_xu_n(t,x)\cdot v
$$
in such a way that 
$$
\int_s^t h_n(\sigma,\gamma(\sigma),\dot \gamma(\sigma))d\sigma
\geq 0
$$
for all absolutely continuous curves $\gamma$ and all $s\leq t$ in $\Rm$.
We deduce that 
 $h_n$ is a non-negative function, and then 
$
 \int h_n d\eta \geq 0.
$
We have the equality (\ref{transeq2}) for $u_n$:
$$
\int_{I\times \Rm^{2d}} \partial_tu_n(t,x)+\partial_xu_n(t,x)\cdot v
d\eta(t,x,v)=0
$$
wich implies that $\int f_nd\eta=\int h_nd\eta\geq 0$.
At the limit $n\lto \infty$, we conclude that $\int fd\eta\geq 0$.
This ends the proof of  Theorem \ref{Young}.\qed

\subsection{Application to the continuity equation}
Young's superposition principle implies elegant results of Ambrosio,
Gigli and Savar\'e  concerning the continuity equation
\begin{equation}\tag{PDE}\label{PDE}
\partial_t \mu + \text{div}(V\mu)=0.
\end{equation}
It is well-known that close relations exist
between (PDE) and the following
\begin{equation}\tag{ODE}\label{PDE}
\dot \gamma(t)=V(t,\gamma(t)).
\end{equation}
More precisely, if $\gamma(t)$ is an absolutely
continuous solution of (ODE), then $\mu:=dt\otimes \delta_{\gamma(t)}$
is a weak solution of (PDE). We call elementary these
solutions.
The relations between (PDE) and (ODE) are enlightened by the following
result, which was obtained by Ambrosio, Gigli and Savar\'e
\cite{AGS,Am:cetraro}, in
the line of anterior works of Smirnov \cite{Si:94} and Bangert
\cite{Ba:99}:

\begin{thm}[Ambrosio, Gigli, Savar\'e]
Let $V:I\times M\lto TM$ be a Borel
time-dependant  vectorfield.
Every probabilty measure $\mu$
on $I\times M$ which solves (PDE)
in the sense of distributions and satisfies
the integrability condition
$$
\int \|V(t,x)\|_x d\mu<\infty
$$ 
is a superposition of 
elementary solutions. 
More precisely, there exists a Borel probability measure $\nu$ on
$\mG(I,M)$ such that $\mu= dt\otimes (ev_t)_{\sharp}\nu$, 
and $\nu$-almost every generalized curve is a
curve and is a   solution of (ODE).
\end{thm}

\proof
In order to see the relation between this result and Young
superposition principle, observe that weak solutions
of (PDE) are in bijection with transport measures which are concentrated
on the graph of $V$ (which is a Borel subset of $X$).
More precisely, if $\eta$ is such a transport  measure, then 
its marginal $\mu$ on $[0,1]\times \Rm^d$
is a weak solution of (PDE).
Conversely, if $\mu$ is a solution of (PDE),
then its lifting to the graph of $V$ is a transport measure.
Now the transport measure $\eta$ associated to the 
solution $\mu$ can be written as a superposition of generalized curves
which are concentrated on the graph of $V$.
But it is obvious that a generalized curve which is concentrated on
the graph of $V$ is nothing but an absolutely continuous solution
of (ODE).

\qed
   
Note that the result can be applied in $\Rm^d$
endowed with the complete metric 
$$
g_x(v,w)=\frac{\langle v,w\rangle}{(1+|x|)^2}.
$$
The integrability condition then reads 
$$
\int_{I\times \Rm^d} \frac{|V(t,x)|}{1+|x|}d\mu(t,x)<\infty
$$
as in \cite{Am:cetraro}.

\subsection{Application to optimal transport}
Let $L$ be a normal integrand.
A generalized curve $\Gamma$ is called minimizing
if it is minimizing the action with fixed boundary points.
If $\eta$ is minimizing the action in 
$\mT_{\mu_i}^{\mu_f}(I,M)$, then $\eta$ can be decomposed into
minimizing generalized curves.
The  decompositions $\nu$ of $\eta$ are  minimizing the action
$$
\int_{\mG(I,M)} \int_{I\times TM} Ld\Gamma d\nu(\gamma)
$$
on the set of probability measures 
$\nu$ on $\mP_1(M)$ such that $(ev_a)_{\sharp} \nu=\mu_i$
and $(ev_b)_{\sharp}\nu=\mu_f$.

If in addition the integrand $L$ is fiberwise convex,
and if there exists a minimizing transport measure
$\eta$ in $\mT_{\mu_i}^{\mu_f}(I,M)$, then there exists
a minimizing transport measure in $\mT_{\mu_i}^{\mu_f}(I,M)$
which is concentrated on the graph of a Borel vector-field $V(t,x)$.
This minimizing measure can be decomposed into minimizing
curves which are  solutions of (ODE).

\section{Holonomic and closed measures}\label{s6}
In the theory of Mather minimizing measures, several spaces of
measures were introduced on $\Tm\times TM$.
In order to be coherent with the exposition of the rest of the present
paper, we shall view them, in an equivalent way, as 
tranport measures in $\mT([0,1],M)$.

\subsection{Closed measures}
They have been used in the context of Lagrangian
dynamics by Bangert in \cite{Ba:99}.
\begin{defn}
A measure $\eta\in \mT([0,1],M)$ is called closed if there exists 
a probability measure $\mu$ on $M$ such that 
$\eta \in \mT_{\mu}^{\mu}([0,1],M)$. 
We denote by $\mF(M)$ the set of closed measures, so that 
$$
\mF(M)=\bigcup_{\mu\in \mP(M)} \mT_{\mu}^{\mu}([0,1],M)
\subset \mT([0,1],M)
$$
\end{defn}

We now expose a superposition principle for closed measures
in the spirit of Smirnov \cite{Si:94}, Bangert \cite{Ba:99}
and  
De Pascal, Gelli and Granieri \cite{DeGeGr}.
Let us first define the set $\mG(\Rm,M)$ 
of measures $\Gamma$
on $\Rm\times TM$ such that, for each $[a,b]\subset \Rm$,
the rescaled restriction
$$
\Gamma_{[a,b]}:=\Gamma_{|[a,b]\times TM}/(b-a)
$$
is a generalized curve in $\mG([a,b],M)$.
Denoting by $d^k$ the distance on $\mG([-k,k],M)$,
we have a distance
$$
d(\Gamma,\Gamma')=\sum_{k=1}^{\infty} 
\frac{d^k
\big(\Gamma_{[-k,k]},\Gamma'_{[-k,k]}\big)}
{2^k}
$$
on $\mG(\Rm,M)$.
Clearly, a sequence $\Gamma^n$ of elemets of $\mG(\Rm,M)$
is converging to $\Gamma$ if and only if we have 
$\Gamma^n_{[a,b]}\lto\Gamma_{[a,b]}$
for each $[a,b]\subset \Rm$.
It is not hard to check that $\mG(\Rm,M)$
is a complete and separable metric space.
Let $\tau:\Rm\times TM\lto \Rm\times TM$
be the translation $(t,x,v)\lmto (t+1,x,v)$.
The map $\tau_{\sharp}: \mG(\Rm,M)\lto \mG(\Rm,M)$
is continuous.
Consequently, the map 
$$
\tau_{\sharp\sharp}: \mP(\mG(\Rm,M))
\lto \mP(\mG(\Rm,M))
$$
is continuous.
A probablity measure $\nu$ on $\mG(\Rm,M)$
is called translation invariant if $\tau_{\sharp\sharp}\nu=\nu$.

For each compact time interval $I$, we denote by 
$$P_I:\mG(\Rm,M)\lto \mG(I,M)$$
the map $\Gamma\lmto \Gamma_I$.
The Borel $\sigma$-algebra
of $\mG(\Rm,M)$ is also the $\sigma$-algebra induced 
by the projections $P_I, I\subset \Rm$.

\begin{thm}
If $\eta$ is a closed measure on $M$, then there exists a 
translation-invariant probabilty measure $\vartheta$
on $\mG(\Rm,M)$ such that 
$$
\int_{[0,1]\times TM} f(t,x,v) d\eta(t,x,v)
=
\int_{\mG(\Rm,M)}\int_{[0,1]\times TM} f(t,x,v) d\Gamma(t,x,v)
d\vartheta(\Gamma)
$$
for each function $f\in L^1(\eta)$.
We call $\vartheta$ a solenoidal decomposition of $\eta$.
\end{thm}

\proof
The proof is based 
on Young's superposition principle and 
on general  constructions of measure theory.
We have $\eta\in \mT_{\mu}^{\mu}$ for some 
probability measure $\mu$ on $M$.
Let $\nu$ be a decomposition of $\eta$
in the sense of Theorem \ref{decomposition}.
We claim that, for each $k\in \Nm$,
there exists a Borel probability measure $\nu^k$ on 
$\mG([0,k],M)$ such that 
$P_{[l,l+1]\sharp \nu^k}=\tau^l_{\sharp} \nu$
for each  $l\in \{0,1,\ldots, k-1\}$.
Then, by standard extension theorems,
(for example Theorem V.4.1 of \cite{Pa})
there exists a unique  probability
$\vartheta$ on $\mG(\Rm,M)$ such that 
$P_{[0,1]\sharp}\vartheta=\nu$, and it is translation invariant.

We have to prove the existence of the measures $\nu^k$.
Let $\nu_x$ be the disintegration of $\nu$
with respect to the map $ev_0$.
In other words, $M\ni x\lmto \nu_x$ is a measurable family
of Borel probabillity measures on $\mG([0,1],M)$ 
such that $\nu_x$ is concentrated on the set of generalized
curves $\Gamma$ which satisfy $ev_0(\Gamma)=x$ and such that
$$
\int_{\mG([0,1],M)} f(\Gamma) d\nu=
\int _M\int_{\mG([0,1],M)}  f(\Gamma) d\nu_x(\Gamma) d\mu(x)
$$
for each bounded Borel function on $\mG([0,1],M)$.
Let us denote 
$$
\mY_1([0,k-1],M)\times \mY_1([k-1,k],M)\lto \mY_1([0,k],M)
$$
$$
(Y,Z)\lmto Y\star Z
$$
the natural gluing.
Note that this map is continuous.
We can now define the  sequence $\nu^k$ by recurrence
setting  $\nu^1:= \nu$
and
$$
\int_{\mG([0,k],M)} f(\Gamma)d\nu^k(\Gamma)=$$
$$
\int_{\mG([0,k-1],M)} 
\int_{\mG([0,1],M)} f(Y\star\tau^{k-1}_{\sharp} Z)
\,d\nu_{ev_{k-1}(Y)}(Z)
  d\nu^{k-1}(Y)
$$
for each bounded continuous function $f$ on $\mG([0,k],M)$.
Note in this expression that $Y\star\tau^{k-1}_{\sharp} Z$
is indeed a generalized curve for $\nu_{ev_{k-1}(Y)}$-almost all
$Z$ because the measure $\nu_{ev_{k-1}(Y)}$ is supported 
on the set of generalized curves $\Gamma$ which satisfy 
$ev_0(\Gamma)=ev_{k-1}(Y)$.
\qed

Let $L:[0,1[\times TM\lto \Rm \cup \{\infty\}$
be a normal integrand. We  extend $L$ by periodicity
to a function on $\Rm\times TM$.
We say that the generalized curve $\Gamma\in \mG(\Rm,M)$
is globally minimizing the action if $\Gamma_{I}$
is minimizing in $\mG(I,M)$
(with fixed endpoints) for each compact inteval $I$.
Similarly, an absolutely continuous curve
 $\gamma:\Rm\lto M$
 is called globally 
  minimizing if it is minimizing the action with fixed
 endpoints on each compact inteval of time.
 If $\eta$ is a closed measure which minimizes the action in $\mF$,
 and if $\vartheta$ is a 
 solenoidal decomposition of $\eta$, then 
 $\vartheta$-almost every generalized curve is globally minimizing.
 If, in addition, the integrand $L$ is fiberwise stricly convex,
 then each minimizing closed measure $\eta$ is concentrated on the 
 graph of a Borel vectorfield $V(t,x)$, this was observed
 in \cite{DeGeGr} and can be proved as the similar statements
 in Section \ref{tonelli}.
 In addition,
 if $\vartheta$ is a solenoidal  decomposition of $\eta$, then 
 $\vartheta$-almost every generalized curve $\Gamma\in \mG(\Rm,M)$
 is a curve, is a solution of (ODE)
 (with the vectorfield $V$ extended to $\Rm\times M$
 by periodicity), and is globally minimizing.
 This property is the generalization in our setting of the theorems
 of Ma\~n\'e \cite{Ma:96}, Bangert \cite{Ba:99}, Fathi and Siconolfi
 \cite{FS} stating,
 under additional assumptions on $L$,
  that minimizing closed measures are invariant.

\subsection{Holonomic measures}
Our analysis of Closed measures makes it well suited
to minimization problems. However,  Ma\~n\'e 
first introduced  in \cite{Ma:96}
the \textit{a priori} smaller set of  holonomic measures.
For historical reasons, 
 we believe it is worth proving
here the equality between holonomic measures and closed measures.
Let $T\in \Nm$ and $\gamma:\Rm\lto M$ be a $T$-periodic absolutely continuous curve. We denote by $\tilde \gamma$
the closed measure defined by
$$
\int_{[0,1]\times TM} fd\tilde \gamma=\frac{1}{T}
\int_0^T f(t-[t],\gamma(t),\dot\gamma(t))dt
$$
where $[t]$ is the integral part of $t$.

\begin{defn}
The set $\mH(M)$ of holonomic measures is the closure, in 
$\mT([0,1],M)$,
of the set of all measures of the form $\tilde \gamma$,
form smooth $T$-periodic curves $\gamma$, $T\in \Nm$.
\end{defn}

\begin{lem}
The set $\mH(M)$ of holonomic measures
can equivalently be defined as the closure in 
$\mT([0,1],M)$,
of the set of all measures of the form $\tilde \gamma$,
for all absolutely continuous $T$-periodic 
curves $\gamma$, $T\in \Nm$.
\end{lem}

\proof
It is sufficient to prove that each measure
$\tilde \gamma$, where $\gamma$ is a $T$-periodic
absolutely continuous curve belongs to $\mH(M)$.
Let $\gamma$ be such a curve.
Let $\gamma_n$ be a sequence of smooth
$T$-periodic curves which converge
to $\gamma$ in $W^{1,1}([0,T],M)$.
Then we prove as in Lemma \ref{continu}
that $\tilde \gamma_n\lto \tilde \gamma$.
\qed

We recall a first remark of Ricardo Ma\~n\'e:

\begin{lem}\label{convex}
The set $\mH(M)$  is convex if $M$ is connected.
\end{lem}

\proof
Let $\eta_1$ and $\eta_2$
be holonomic measures, and let $\lambda_1$ and $\lambda_2$
in $[0,1]$ be such that $\lambda_1+\lambda_2=1$.
We want to prove that $\lambda_1\eta_1+\lambda_2\eta_2$
is holonomic.
Since $\eta_i$ is holonomic,
there exists a sequences of integers $T_n^i$ and a sequences
of smooth  curves $\gamma_n^i(t):\Rm\lto M$
of period $T_n^i$ such that $\tilde \gamma_n^i\lto \eta_i$.
By possibly replacing the periods $T^i_n$ by multiples,
we can suppose without loss of generality that 
$$
T_n^1/T_n^2\lto \lambda_1/\lambda_2.
$$
Let $(0,x^i,v^i)$ be a point in the support of $\eta^i$.
Since $\tilde \gamma^i _n\lto \eta^i$,
 there exists a sequence $t^i_n$ of times such that 
$$(t_n^i -[t_n^i],\gamma(t_n ^i), \dot \gamma(t_n^i))
\lto (0,x^i,v^i).
$$ 
We can suppose   that 
$t^i_n\lto 0$ by replacing the curves 
$ \gamma^i_n(t)$ by $\gamma^i_n(t-[t^i_n])$.
We consider  the sequence $\gamma_n$ of absolutely continuous 
curves of period $T_n^1+2+T_n^2$ such that 
$\gamma_n=\gamma_n^1$ on $[t_n^1,T_n^1+t_n^1]$, 
$\gamma_n=\gamma_n^2$ on $[1+T_n^1+t_n^2,1+T_n^1+T_n^2+t_n^2]$
and $\gamma_n$ is a minimizing geodesic
on the remaining intervals.
It is not hard to see that
 $$\tilde \gamma_n \lto \lambda_1\eta_1+\lambda_2\eta_2$$
 as $n\lto \infty$,
so that this measure is holonomic.
\qed

The following result
 is a piece of unproved folklore:

\begin{thm}\label{closed}
If $M$ is a compact connected manifold, then each 
closed measure is Holonomic :
$$
\mH(M)=\mF(M).
$$
\end{thm}

\proof
By Lemma \ref{convex} and Proposition \ref{dual}, it is enough
to prove that if $f\in C_1([0,1]\times TM)$ satisfies 
$$
\int_0^T f(t-[t], \gamma(t),\dot \gamma(t))dt
\geq 0
$$
for each $T\in \Nm$ and each absolutely continuous 
$T$-periodic curve $\gamma:\Rm\lto M$,
then $\int fd\eta\geq 0$ for each closed measure $\eta$.
We fix an integrand $f\in C_1([0,1]\times TM)$, and extend
 $f_{|[0,1[}$
to a $1$-periodic function on $\Rm\times TM$
without changing the name.
We have 
$$
\int _0^T f(t, \gamma(t),\dot \gamma(t))dt
\geq 0
$$
for each $T\in \Nm$ and each absolutely continuous 
$T$-periodic curve $\gamma$.
As in the proof of Young's principle, we consider the value function
$u:[0,\infty)\times M\lto \Rm$
defined by
$$
u(t,x):=\inf_{\gamma(t)=x}
 \int_0^t f(s,\gamma(s),\dot\gamma(s)) ds
$$
where the infimum is taken on the set of absolutely continuous
curves $\gamma:[0,t]\lto M$ which satisfy $\gamma(t)=x$.

\begin{lem}
The function $u$ is upper  semi-continuous and locally bounded on 
$[0,\infty)\times N$. 
In addition, it is bounded from below.
\end{lem} 

\proof
The proof of the first part is similar to the proof of Lemma 
\ref{usc}.
We prove that the value function is bounded from below.
There exists a constant $C$ such that
$$
\int _0^s f(t,\gamma(t),\dot \gamma(t)) dt
\geq -C
$$
for each $s\geq 0$ and each absolutely continuous
curve $\gamma:[0,s]\lto M$.
Indeed, we can consider the periodic curve
$x(t):[0,[s]+2]\lto M$ such that
$x=\gamma$ on $[0,s]$
and 
$
x(t)
$
is a minizing geodesic on
  $[s,[s]+2]$
  between $\gamma(s)$ and $\gamma(0)$.
We have
$$
0\leq \int_0^{[s]+2}f(t,x(t),\dot x(t))dt
\leq \int_0^s f(t,\gamma(t),\dot \gamma(t)) dt
+2\|f\|_1\big(1+D(\gamma(0), \gamma(s))\big),
$$
where $D$ is the Riemannian distance on $M$, 
which is bounded.
From the definition of $u$, it follows that $u(s,x)\geq -C$
for each $(s,x)$ in $[0,\infty)\times M$.
\qed

As in the proof of Young's principle, we have:

\begin{lem}
There exist sequences $u_n:[1,\infty) \times M\lto \Rm$ and 
$f_n:\Rm\times TM\lto \Rm$ of functions such that:
\begin{itemize}
\item The functions $f_n$ are $1$-periodic in $t$.
They are countinuous and satisfy a uniform estimate
$|f_n(t,x,v)|\leq C(1+\|v\|_x)$. Finally, we have $f_n\lto f$
almost everywhere.
\item The functions $u_n$ are smooth, locally bounded
 and bounded from below.
\item
 The inequality
$$
 u_n(t,\gamma(t))-u_n(s,\gamma(s))\leq
\int_s^t f_n(\sigma,\gamma(\sigma),\dot\gamma(\sigma))d\sigma
$$
holds for each $1\leq s\leq t$ in $\Rm$ and each absolutely continuous curve
$\gamma:\Rm\lto M$.
\end{itemize}
\end{lem}

\proof
We regularize as in the proof of Young's principle.
There is however a small difficulty related to the fact that 
we now work on a manifold.
In order to solve this difficulty, we embed $M$ as a Riemannian
submanifold of some Euclidean space $\Rm^d$ -- one could
also regularize in a more intrinsic way in the spirit of De Rham
\cite{DR}.
Then, we consider a  tubular neighborhood 
$U$ of $M$ in $\Rm^d$ and the associated smooth projection 
$\pi:U\lto M$.
We set, for $t\geq 1$ and $(x,v)\in TM$,
$$
u_n(t,x)
=
\int_{\Rm\times \Rm^d}
u(t-\sigma,\pi(x-y)) \rho_n(\sigma,y)d\sigma dy
$$
and
$$
f_n(t,x,v)=
\int_{\Rm\times \Rm^d}
f(t-\sigma,\pi(x-y),d\pi_{(x-y)}\cdot v)\rho_n(\sigma,y)d\sigma dy.
$$
Let $\gamma :\Rm\lto M$ be an absolutely continuous curve.
For each fixed small $y\in \Rm^d$,
the curve $\pi(\gamma(t)-y)$ is absolutely continuous, and we have
$$
u(t-\sigma, \pi(\gamma(t-\sigma)-y))
-u(s-\sigma, \pi(\gamma(s-\sigma)-y))$$
$$
\leq 
\int_{s}^t
 f\Big(\zeta-\sigma, \pi(\gamma(\zeta-\sigma)-y),
d\pi_{(\gamma(\zeta-\sigma)-y)}\cdot \dot \gamma(\zeta-\sigma)\Big)
d\zeta
$$
for all small $s$.
The third point of the Lemma follows by integration.
The other points are standard.
\qed

Let $\eta\in \mT_{\mu}^{\mu}([0,1],M)$ be a closed measure.
We want to prove that $\int fd\eta\geq 0$.
We see as in the proof of Young's principle
that
$$
\int _M u_n(i+1,x)-u_n(i,x)d\mu(x)\leq
\int_{[0,1]\times TM} f_n d\eta
$$
for each integer $i\geq 1$.
By summation, we obtain, for each $T\in \Nm$,
$$
\int _M\frac{u_n(T+1,x)-u_n(1,x)}{T}d\mu(x)\leq
\int_{[0,1]\times TM} f_n d\eta.
$$
At the limit $T\lto \infty$, we obtain that 
$\int f_nd\eta\geq 0$, and then at the limit $n\lto \infty$,
we get $\int fd\eta\geq 0$, as desired.
This ends the proof of Theorem \ref{closed}.
\qed


\begin{appendix}

\section{Kantorovich-Rubinstein space}
 Good references for the material exposed here are \cite{AGS}
 and \cite{Vi:03}.
 Let $(X,d)$ be a complete and separable metric space.
Let $\mP_1(X)$ be the set of Borel probability measures on $X$ with 
finite first moment, that is the set of Borel probability
measures $\mu$ on $X$   such that the integral
$$
\int_X d(x_0,x) d\mu(x)
$$
is finite for one (and then each) point $x_0\in X$.

A coupling between two probability measures $\mu$ and 
$\eta$ is a probability measure $\lambda$ on $X^2$
whose marginals are $\mu$ and $\eta$, or in other words such that
$$
\int_{X^2}f(x)+g(y) d\lambda(x,y)
=
\int_X f(x)d\mu(x)+\int_X g(y) d\eta(y)
$$
for all continuous functions $f$ and $g$ on $X$.

We recall the definition of the Kantorovich-Rubinstein distance
$d$ on $\mP_1(X)$:
$$
d_1(\mu,\eta)=\min _{\lambda} \int_{X\times X} d(x,y)
 d\lambda(x,y)
$$
where the minimum is taken on the set of couplings
$\lambda$ between $\mu$ and $\eta$.

Let us denote by $C_1(X)$ the set of continuous functions $f$  on $X$
such that 
$$
\sup_{x\in X} \frac{|f(x)|}{1+d(x_0,x)}<\infty
$$
for one (and then any) point $x_0\in X$.
The topology on $\mP(X)$ defined by 
the distance $d$
is precisely the weak topology associated to the linear forms
$\mu\lmto \int fd\mu$, $f\in C_1(X)$.
In other words, we have $d(\mu_n,\mu)\lto 0$
if and only if
$$
\int fd\mu_n\lto \int fd\mu
$$
for all $f\in C_1(X)$.
There is an interesting duality formula for the distance:
$$
d_1(\mu,\eta)=\sup _{f} \int_{X} f(x)d(\mu-\eta)(x)
 $$
where the supremum is taken on the set of $1$-Lipschitz functions
$f:X\lmto \Rm$. An important remark is that, if the distance $d$ 
on $X$ is bounded, then the associated Kantorovich-Rubinstein 
space is just the space $\mP(X)$ of all Borel probability 
measures on $X$ endowed with the narrow topology.
Since it is always possible to replace a given distance $d$
by another distance which is bounded and generates the same topology,
our discussion includes the study of the narrow topology
on $\mP(X)$.

The metric space  $(\mP_1(X),d)$
 is complete and separable, see \cite{AGS}.
The relatively compact subsets of $\mP_1(X)$ are those which are 
$1$-tight:
 
 \begin{defn}
 The subset $Y\subset \mP(X)$
 is called $1$-tight if one of the following equivalent properties holds:
 \begin{itemize}
 \item For each $\epsilon>0$, there exists a compact set
 $K\subset X$ and a point $x_0$ such that
 $\int_{X-K} (1+d(x_0,x))d\mu \leq \epsilon$
 for each $\mu \in Y$.
 \item There exists a  function $f:X\lto [0,\infty]$
 whose sublevels are compact,
 a constant $C$ and a point $x_0$ such that 
 $\int_X (1+d(x_0,x))f(x)d\mu\leq C$ for each $\mu \in Y$.
\item The family $Y$ is tight with uniformly integrable first 
moment.
The first means that, for each $\epsilon>0$, there exists a compact set
 $K\subset X$  such that
 $\mu (X-K) \leq \epsilon$
 for each $\mu \in Y$.
 The second means that for each $\epsilon>0$, there exists
 a ball $B$ in $X$ such that 
 $\int_{X-B} d(x_0,x)d\mu \leq \epsilon$
 for each $\mu \in Y$.
 \end{itemize}
\end{defn}
 Note that $1$-tightness is just tightness if the distance $d$ is bounded.

\begin{lem}
A sequence $\mu_n$ converges to $\mu$ in $\mP_1(X)$ 
if and only if the family $\{\mu_n,n\in \Nm\}$
is $1$-tight and if 
$\mu_n$ narrowly converges to $\mu$, which means that
$\int fd\mu_n\lto \int fd\mu$
for each bounded  continuous function $f$.
It is enough that the family $\mu_n$ is converging narrowly
to $\mu$ and has uniformly integrable first moment. 
\end{lem}

Let us now assume that $X$ is a finite dimensional manifold.
\begin{lem}
A sequence $\mu_n$ converges to $\mu$ in $\mP_1(X)$ 
if and only if the sequence has uniformly integrable first moment
and converges to $\mu$ in the sense of distributions.
\end{lem}

Still assuming that $X$ is a manifold, we finish with the following:

\begin{lem}
Let $\mu_t,t\in I$ be a measurable family of probability measures
on $X$, where $I$ is an interval of $\Rm$.
In order that $\mu_t$ is equal almost everywhere to a narrowly 
continuous map, it is enough that, for each compactly supported
smooth funtion $f:X\lmto \Rm$, the function $t\lmto \int fd\mu_t$
is equal almost everywhere to a continuous function.
\end{lem}

\section{Superpositions}
We continue with the notations of the first appendix.
Let $\nu$ be a Borel probability measure on the 
complete metric space $\mP_1(X)$.
We say that $\nu$ represents the measure 
$\eta\in \mP_1(X)$ if the equality
\begin{equation}\label{superposition}
\int _X fd\eta=
\int_{\mP_1(X)}\int_X  fd\mu d\nu(\mu)
\end{equation}
holds for each function $f\in L^1(\eta)$.
Let us first check that the right hand side is meaningful:

\begin{lem}
The linear map $\mu \lmto \int_X  fd\mu$ is Borel measurable
on $\mP_1(X)$ when $f$ is a non-negative Borel function on $X$.
Each probability measure  $\nu$ on $\mP_1(X)$
represents one (and only one) element $\eta\in \mP(X)$.
We have $\eta\in \mP_1(X)$ if and only if
$$\int_{\mP_1(X)}\int_X  d(x_0,x)d\mu(x) d\nu(\mu)<\infty
$$
for one point $x_0$.
\end{lem}

\proof
It is clear that (\ref{superposition})
defines a (unique) Borel measure $\eta$ if it  is meaningful 
for each non-negative Borel function.
So we have to prove the first statement. 
Since the conclusion holds when $f$ is continuous and bounded
(for then the map $\mu\lmto \int fd\mu$ is continous),
it  is a consequence of the following standard statement.
\qed

\begin{lem}
Let $E$ be a vector space of real-valued functions on $X$.
Assume that $E$ contains all  bounded continuous functions and is
closed under monotone convergence. Then $E$ contains all
non-negative Borel functions.
\end{lem}

\proof
Let $\mB$ be the set of subsets of $X$ whose caracteristic function
belongs to $E$. It is not hard to see that $\mB$ contains
closed sets, that it is closed under increasing union,
and that if $A\subset B$ are two elements of $\mB$, then 
$B-A$ is an element of $\mB$.
The classical Dynkin class theorem then implies that 
$\mB$ contains all the Borel sets.
But then $E$ contains all Borel non-negative functions.
\qed

This statement also implies:

\begin{lem}
In order that (\ref{superposition}) holds for each function 
$f\in L^1(\eta)$, it is sufficient that it holds for all 
bounded continuous functions.
\end{lem}

\begin{prop}
Let $\mG$ be a closed subset of $\mP_1(X)$,
and let $\mT$ be the closed convex envelop of $\mG$ in $\mP_1(X)$.
Each measure $\eta\in \mT$ is represented by a measure
$\nu$ which is supported on $\mG$
(we say that $\mu$ is a superposition of elements of $\mG$). 
\end{prop}

\proof
Let us consider the set $\mS$ of elements of $\mP(X)$ which
are superpositions of elements of $\mG$.
It is obvious that the set $\mS$ is convex, and contains
$\mG$.
So we have to prove that this set is closed.
Let us consider a sequence $\eta_n$ in $\mS$,
which has a limit $\eta$ in $\mP_1(X)$.
There exists a sequence $\nu_n$ of Borel probability
measures on $\mP_1(X)$ which represents $\eta_n$.
Since the family $\{\eta,\eta_1, \ldots, \eta_n,\ldots\}$
is compact in $\mP_1(X)$, it is $1$-tight, hence there exists
a function $f:X\lto [0,\infty]$
whose sublevels $f^{-1}([0,c])$ are compact and such that 
the integral
$$
\int_X (1+d(x_0,x))f(x)d\eta_n(x)=
\int_{\mP_1(X)} \int_{X} (1+d(x_0,x)) f(x)d\mu(x) d\nu_n(\mu)
$$
is a bounded sequence.
The map
$$
\mu\lmto \int_{X} (1+d(x_0,x)) f(x)d\mu(x) 
$$
has compact sublevels on $\mP_1(X)$, hence the  boundedness
of the sequence above  implies
that the sequence $\nu_n$ is a tight sequence
of probability measures on $\mP_1(X)$.
By the standard Prohorov theorem, we can assume that 
$\nu_n$ has a limit $\nu$ for the narrow topology, 
which means that
$$
\int_{\mP(X)}F(\eta) d\nu_n(\eta)
\lto \int_{\mP(X)}F(\eta) d\nu(\eta)
$$
for each bounded and continuous function
$F$ on $\mP_1(X) $.
For each continuous and bounded function $f$ on $X$,
the affine  function $\mu \lmto \int_Xf d\mu$ is continuous
and bounded on $\mP(X)$,
hence
$$
\int_{\mP(X)}\int _X fd\mu  d\nu_n(\mu)
\lto \int_{\mP(X)}\int _X fd\mu d\nu(\mu).
$$
Recalling that 
$$
\int_{\mP(X)}\int _X fd\mu  d\nu_n(\mu)=
\int _X fd\eta_n\lto \int _X fd\eta,
$$
we conclude that
$$
 \int _X fd\eta=
 \int_{\mP(X)}\int _X fd\mu d\nu(\mu)
$$
for each bounded and continuous function $f$ on $X$.
This implies that  $\nu$ represents 
$\mu$. Since the measures $\nu_n$ are supported on
the closed ste  $\mG$, the limit
$\nu$ is supported on $\mG$.
We have proved that $\mu\in \mS$.
\qed

We  finish with  an obvious remark on closed convex subsets of
$\mP_1(X)$.

\begin{prop}\label{dual}
Let $\mC$ be a closed convex subset of $\mP_1(X)$, and
let $\mC^+$ be the set of functions $f\in C_1(X)$ such that
$\int_X fd\mu\geq 0$ for each $\mu\in \mC$.
Then $\mC$ is the set of measures $\mu\in \mP_1(X)$
such that
$\int_X fd\mu\geq 0$
for each $f\in \mC^+$.
\end{prop}

\end{appendix}


\begin{thebibliography}{99}


\bibitem{Am:BV}
L. Ambrosio, 
\textit{Lecture notes on transport equation and Cauchy
problem for BV vector fields and applications.}

\bibitem{Am:cetraro}
L. Ambrosio, 
\textit{Transport equation and Cauchy problem for non-smooth
vector fields}


\bibitem{AGS}
L. Ambrosio,  N. Gigli and  G. Savar\'e, 
Gradient flows, \textit{Lectures in Math. ETH Z\" urich}, 
Birkh\"auser (2005).

\bibitem{Ba:99}
V. Bangert, 
\textit{Minimal measures and minimizing closed normal one-currents},
GAFA, \textbf{9} (1999), 413-427.

\bibitem{BB}
J. D. Benamou, Y.  Brenier, 
\textit{A computational fluid mechanics solution to the Monge-Kantorovich mass transfer problem.} Numer. Math. 
\textbf{84} (2000), no. 3, 375--393. 

\bibitem{BeLa:73}
H. Berliocchi  and J.-M. Lasry,
\textit{Int\'egrandes normales et mesures param\'etr\'ees en calcul des variations,}
Bull. Soc. Math. France \textbf{101} (1973), 129--184.

\bibitem{BeBu}
P. Bernard and B. Buffoni,
\textit{Optimal mass transportation and Mather theory}, 
JEMS \textbf{9} (2007), no. 1,  85--121.


\bibitem{Br:89}
Y. Brenier, \textit{The least action principle and 
the related concept of generalized flows for incompressible perfect fluids,}
J. A. M. S. \textbf{2} (1989), 225--255.

\bibitem{Br:93}
Y. Brenier, \textit{The dual least action problem for an ideal,
incompressible fluid}
Arch. Rat. Mech. Anal. \textbf{122} (1993), 323--351.

\bibitem{Br:99}
Y. Brenier, \textit{Minimal geodesics on groups of volume-preserving
maps and generalized solutions of the Euler equations,}
Comm. Pure. Ap. Math. \textbf{52} (1999), 411--452. 

\bibitem{Br:03}
 Y. Brenier, \textit{Extended Monge-Kantorovich theory.}
 Optimal transportation and applications (Martina Franca, 2001), 91--121, Lecture Notes in Math., 
 \textbf{1813} , Springer, Berlin, 2003.



\bibitem{DeGeGr} L. De Pascale, M. S. Gelli  and L. Granieri,
\textit{ Minimal measures, one-dimensional currents and the
Monge-Kantorovich problem},
Calc. Var.  Part. Diff. Eqn.
\textbf{27} (2006), no. 1, 1--23.


\bibitem{DR}
G. De Rham,
\textit{Vari\'et\'es diff\'erentielles,} Hermann, 1960.

\bibitem{FS}
A. Fathi and A. Siconolfi,
\textit{Existence of $C^1$ critical 
sub-solutions of the Hamilton-Jacobi equation,}
 Invent. Math. \textbf{155}
   (2004), no. 2, 363--388.


\bibitem{Ma:91}
J. N. Mather, \textit{Action minimizing invariant measures for positive
definite Lagrangian systems}, Mathematische Zeitschrift,
Math. Z. {\bf 207} (1991), 169-207.

\bibitem{Ma:96}
  R. Ma\~n{\'e}, \textit{Generic properties and problems of minimizing measures of
   {L}agrangian systems}, Nonlinearity \textbf{9} (1996), no.~2, 273--310.
\bibitem{Pa}
K. R. Parthasarathy,
Probability measures on metric spaces,
Academic Press (1967).

\bibitem{Si:94} S.K. Smirnov,\textit{
Decomposition  of solenoidal vector charges into elementary 
solenoids and the structure of normal one-currents}
St. Petersbourg Math. J. \textbf{5} (1994), 841-867,
\bibitem{Va:94}
M. Valadier,
\textit{A course on Young measures,}
Rend. Istit. Mat. Univ. Trieste
\textbf{24} (1994), 349--394.


\bibitem{Vi:03}
C. Villani,
Topics in optimal transportation, American Mathematical Society,
Providence, Rhode Island, 2003.

\bibitem{Yo}
L. C. Young,
Lectures on the calculus of variations and optimal control 
theory, second edition, Chelsea (1980).


\end{thebibliography}
\end{document}